\documentclass[english]{amsart}

\usepackage{latexsym}

\usepackage{amsmath}
\usepackage{rotating}
\usepackage{amsfonts}
\usepackage{amsthm}
\usepackage{amssymb}
\usepackage{verbatim}
\usepackage{psfrag}
\usepackage{graphicx}

\newtheorem{lemm}{Lemma}[section]
\newtheorem{prop}[lemm]{Proposition}
\newtheorem{conj}[lemm]{Conjecture}
\newtheorem{defi}[lemm]{Definition}
\newtheorem{coro}[lemm]{Corollary}
\newtheorem{theo}{Theorem}

  \usepackage[all]{xy}
\xyoption{matrix}

\newcommand{\proofend}{\hfill $\square$}
\newcommand{\car}{\mathrm{char}}

\newcommand{\K}{\mathbb{K}}
\newcommand{\C}{\mathbb{C}}

\newcommand{\KK}{\overline{\mathbb{K}}}

\newcommand{\Aut}{\mathrm{Aut}}

\newcommand{\Pn}{\mathbb{P}^2}

\newcommand{\drawat}[3]{\makebox[0pt][l]{\raisebox{#2}{\hspace*{#1}#3}}}

\title{The correspondence between a plane curve and its complement}
\author{J\'er\'emy Blanc}
\thanks{The author acknowledge support from the Swiss national science foundation}
\address{J\'er\'emy Blanc,
Universit\'e de Grenoble I,
UFR de Math\'ematiques,
UMR 5582 du CNRS, 
Institut Fourier, BP 74, 
38402 Saint-Martin d'H\`eres, France}

\begin{document}
\maketitle

\begin{abstract}
Given two irreducible curves of the plane which have isomorphic complements, it is natural to ask whether there exists an automorphism of the plane that sends one curve on the other. 

This question has a positive answer for a large family of curves and H.~Yoshihara conjectured that it is true in general. We exhibit counterexamples to this conjecture, over any ground field. In some of the cases, the curves are isomorphic and in others not; this provides counterexamples of two different kinds.

Finally, we use our construction to find the existence of surprising non-linear automorphisms of affine surfaces.
\end{abstract}


\section{Introduction}
In this article, $\K$ is any field, and all surfaces are algebraic affine or projective surfaces, defined over $\K$.

\subsection{The conjecture}
To any irreducible curve  $C \subset \Pn=\Pn_{\K}$ we can associate its complement, the affine surface $\Pn\backslash C$ (such affine surfaces have been a subject of research for many years, see \cite{bib:GiDa75}, \cite{bib:Iit77}, \cite{bib:Yos79}, \cite{bib:Miy81}, \cite{bib:Miy01}, \cite{bib:Kis01}, \cite{bib:Koj05}, ...). If two such curves $C,D$ are projectively equivalent -- i.e.~if some automorphism of the projective plane $\Pn$ sends $C$ on $D$ -- then clearly $\Pn\backslash C$ is isomorphic to $\Pn\backslash D$. It is natural to ask whether the converse is true. In $1984$, Hisao Yoshihara made the following conjecture.

\begin{conj}[\cite{bib:Yos84}]\label{Conj:Yoshi}
Let $C\subset \Pn_{\K}$ be an irreducible curve and assume that $\K$ is algebraically closed of characteristic $0$. Suppose that $\Pn\backslash C$ is isomorphic to $\Pn\backslash D$ for some curve $D$. Then $C$ and $D$ are projectively equivalent.
\end{conj}
In \cite{bib:Yos84}, it was proved that the conjecture is true for a large family of curves~$C$. We briefly recall these results in Section~\ref{Case:True}, and extend some of them to any field~$\K$. 
Then, we provide a family of counterexamples to the conjecture, over any field~$\K$, and prove the following result.
\begin{theo}\label{Thm:CountExConj}
For any field~$\K$ with more than two elements, there exist two curves $C,D\subset \Pn_{\K}$, irreducible over the algebraic closure of $\K$, such that the following two assertions are true:
\begin{enumerate}
\item
the affine surfaces
 $\Pn\backslash C$ and $\Pn\backslash D$ are isomorphic;
 \item
 no automorphism of $\Pn$ sends $C$ on $D$.
 \end{enumerate}
Furthermore, there are examples where $C$ and $D$ are isomorphic and examples where they are not. 
\end{theo}

 Observe that Theorem~\ref{Thm:CountExConj} yields the existence of isomorphic affine surfaces having a projective completion in isomorphic projective surfaces by irreducible non-isomorphic curves. Such examples were, as far as we are aware, not known before.
 
Recall that a curve $C$ is of type $I$ if there exists some point $a\in C$ such that $C\backslash a$ is isomorphic to the affine line.
The problem stated above is related to another conjecture, namely:

\begin{conj}[\cite{bib:Yos85}, page 101]\label{Conj:Yoshi2}
If $C\subset \Pn_{\C}$ is an irreducible curve, which is neither of type $\mathrm{I}$ nor a nodal cubic curve, then any automorphism of $\Pn\backslash C$ extends to an automorphism of $\Pn$.
\end{conj}

The construction we provide to prove Theorem~\ref{Thm:CountExConj}  will also provide counterexamples to Conjecture~\ref{Conj:Yoshi2}, extending furthermore the possibilities for the base field.
\begin{theo}\label{Thm:SecondConj}
Assume that the characteristic of $\K$ is not $2$. Then, there exists a curve $C\subset \Pn_{\K}$, irreducible over the algebraic closure of $\K$, of degree $39$, that is not of type $I$, and there exists an automorphism of $\Pn_{\K}\backslash C$ that does not extend to $\Pn_{\K}$.
 \end{theo}

\subsection{The construction}\label{Sec:ConstructionIntro}
Here, we briefly describe our construction, which will be explained more precisely in Section~\ref{Sec:ConstrDetails}. We denote by $\Delta\subset \Pn$ the union of three general lines and choose two quartics $\Gamma_1,\Gamma_2$ that intersect $\Delta$ in a particular manner. We construct a birational morphism $\pi:X\rightarrow \Pn$ that is a sequence of blow-ups of points that belong, as proper or infinitely near points, to  $\Delta\cap\Gamma_1$ or $\Delta\cap\Gamma_2$. Then, we find a reducible curve $R\subset \pi^{-1}(\Delta)$ such that for $i=1,2$, the curve $R\cup{\tilde{\Gamma_i}}$ is contractible via a birational morphism $\eta_i:X\rightarrow \Pn$ (where $\tilde{\Gamma_i}$ is the strict transform of $\Gamma_i$ on $X$). The birational map $\varphi=\eta_1\circ{\eta_2}^{-1}$ restricts to an isomorphism from $\Pn \backslash \eta_2(\tilde{\Gamma_1})$ to  $\Pn \backslash \eta_1(\tilde{\Gamma_2})$.

\begin{equation}\label{Eq:DiagrammFond}
\xymatrix@R=1mm@C=2cm{& &  \Pn\\
\Pn&X\ar[l]_{\pi}\ar@/_/[rd]^{\eta_2}\ar@/^/[ru]^{\eta_1}&\\
&&\Pn\ar@/^1pc/@{-->}[uu]_{\varphi}
}\end{equation}



In our construction, the curves $\Gamma_1$ and $\Gamma_2$ depend on parameters. For general values of these parameters, the curves $\eta_2(\tilde{\Gamma_1})$ and $\eta_1(\tilde{\Gamma_2})$ are not projectively equivalent, which yields the proof of Theorem~\ref{Thm:CountExConj}.
For special values of the parameters, there exists some automorphism $\psi$ of $\Pn$ that sends  $\eta_2(\tilde{\Gamma_1})$ on $\eta_1(\tilde{\Gamma_2})$. Thus, $\varphi\circ \psi^{-1}$ is an automorphism of $\Pn\backslash \eta_1(\tilde{\Gamma_2})$ that does not extend to an automorphism of $\Pn$, which proves Theorem~\ref{Thm:SecondConj}. 
\subsection{Outline of this article}
In Section~\ref{Case:True}, we prove that Conjectures~\ref{Conj:Yoshi} and \ref{Conj:Yoshi2} are true for "most" kinds of curves. In Section~\ref{Sec:ConstrDetails}, we describe precisely the construction announced in (\ref{Sec:ConstructionIntro}). Finally in Section~\ref{Sec:Comparison} we prove that neither of the curves constructed is of type $\mathrm{I}$, and decide when the curves obtained are projectively equivalent or isomorphic, which yields the proofs of Theorems~\ref{Thm:CountExConj} and \ref{Thm:SecondConj}.
\subsection{Aknowledgements}
The author presented the results of this article in Dijon and Gen\`eve. He would like to express his sincere gratitude to the members of these institutes for valuable questions which helped him to improve the exposition of this paper, with special thanks to Adrien Dubouloz and Thierry Vust.
\section{Cases in which the conjectures are true}\label{Case:True}
In this section, we prove that the conjectures are true for most curves, and recall some classical results. We will denote the algebraic closure of $\K$ by $\KK$.
\begin{defi}\label{Defi:TowerResolution}
We say that a birational morphism $\chi:S\rightarrow \Pn$ is a \emph{$n$-tower resolution} of a curve $C\subset \Pn$ if 
\begin{enumerate}
\item
the map $\chi$ decomposes as $\chi=\chi_m\circ\chi_{m-1}\circ... \circ \chi_1$, for some integer $m\geq 0$, where $\chi_i$ is the blow-up of a point $p_i$ and  $\chi_{i-1}(p_i)=p_{i-1}$ for $i=2,..,m$;
\item
the strict transform of the curve $C$ on $S$ is a curve that is smooth, irreducible over $\KK$, isomorphic to $\mathbb{P}^1$, and of self-intersection $n$.
\end{enumerate}
\end{defi}
Note that if a curve admits a $n$-tower resolution, it admits a $m$-tower resolution for any $m\leq n$. Next, we remind the reader of a simple but useful lemma, obvious for the specialist.
\begin{lemm}
Let $C\subset \Pn$ be a curve irreducible over $\KK$, and let $\psi:\Pn\backslash C\rightarrow \Pn\backslash D$ be an isomorphism, where $D$ is some curve of $\Pn$.

Then, either $\psi$ extends to an automorphism of $\Pn$ (and in particular $C$ and $D$ are projectively equivalent), 
or there exist two birational morphisms $\chi,\epsilon:S\rightarrow \Pn$ satisfying the following conditions:
\begin{enumerate}
\item
$\chi$ (respectively $\epsilon$) is a $(-1)$-tower resolution of $C$ (respectively of $D$);
\item
$\chi$ is a minimal resolution of the indeterminacies of $\psi$ and $\psi\chi=\epsilon$.
\end{enumerate}
\end{lemm}
\begin{proof}
In this proof, we consider our algebraic varieties over the field $\KK$, remembering that these are defined over the subfield~$\K$. 
We extend $\psi$ to a birational transformation $\overline{\psi}$ of $\Pn_{\KK}$, which is defined over the field~$\K$. Then, there exists a birational morphism $\overline{\chi}:S_{\KK}\rightarrow \Pn_{\KK}$, also defined over $\K$, that is a minimal resolution of the indeterminacies of $\overline{\psi}$. We denote the birational morphism $\overline{\psi}\circ \overline{\chi}$ by $\overline{\epsilon}$  and denote by $E$ (respectively $F$) the set of irreducible curves of $S_{\KK}$ that are collapsed by $\overline{\chi}$ (respectively by $\overline{\epsilon}$). Since $\psi$ is an isomorphism of $\Pn\backslash C$ to $\Pn\backslash D$, and under the assumption that $\overline{\psi}$ is not an automorphism of $\Pn_{\KK}$, the map $\overline{\psi}$ collapses exactly one irreducible curve of $\Pn_{\KK}$, which is the extension of $C$ as $\overline{C}\subset \Pn_{\KK}$. This means that the set $F\backslash E$ consists of a single element, which is the strict transform of $\overline{C}$; since the sets $E$ and $F$ have the same number of curves, the set $E\backslash F$ also consists of a single element. This element has to be the strict transforms on $S_{\KK}$ of the extension $\overline{D}$ of the curve $D$. 
The resolution of $\overline{\psi}$ by $\overline{\chi}$ and $\overline{\epsilon}$ being minimal, every irreducible curve of $E\cap F$ has self-intersection $\leq -2$; this implies that the strict transforms of $\overline{C}$ and $\overline{D}$ on $S_{\KK}$ are $(-1)$-curves, i.e.~both are smooth, irreducible, isomorphic to $\mathbb{P}^1$ and of self-intersection $-1$.

 The fact that only one irreducible curve collapsed by $\overline{\chi}$ (respectively by $\overline{\epsilon}$) has self-intersection $-1$ implies that $\overline{\chi}$ is a tower resolution of $\overline{C}\subset \Pn_{\KK}$ (respectively of $\overline{D}\subset \Pn_{\KK}$). 
 Since the set of points blown-up by both morphisms is invariant under the action of $\mathrm{Gal}(\KK/\K)$, and since no two points belong to the same surface, each point is defined over $\K$. Consequently, reducing the ground field to $\K$, we find birational morphisms $\chi$ and $\epsilon$ that are tower resolutions of $C$ and $D$ respectively.
\end{proof}
\begin{coro}
Conjectures~\ref{Conj:Yoshi} and \ref{Conj:Yoshi2} are true for any base field~$\K$ and any curve $C\subset \Pn$, irreducible over $\KK$, that does not admit a $(-1)$-tower resolution.

In particular, both conjectures are true if $C$ is not rational or if $C$ has more than two singular points over $\KK$.\proofend
\end{coro}
The conjectures are thus true for a large family of curves. Among curves admitting a tower resolution, curves of type $\mathrm{I}$ or $\mathrm{II}$ are the most natural to deal with. We remind the reader of some results on this subject.

\begin{defi}
A curve $C\subset \Pn$ is of type $\mathrm{I}$ (respectively of type $\mathrm{II}$) if there exists a point $a\in C$ (respectively a line $L\subset \Pn$) such that $C\backslash a$ (respectively $C\backslash L$) is isomorphic to the affine line.
\end{defi}
Any curve of type $\mathrm{II}$ is of type $\mathrm{I}$ and it is difficult (but possible) to find curves of type $\mathrm{I}$ that are not of type $\mathrm{II}$ \cite{bib:Yos83}. A curve is of type $\mathrm{II}$ if and only if it is the image of a line by an automorphism of $\Pn\backslash L$, where $L$ is a line \cite{bib:AbhyMoh}.
 Furthermore, any curve of type $\mathrm{II}$ admits a $n$-tower resolution, for some positive integer $n\geq 3$ \cite{bib:Yos87}. The following result gives another evidence to Conjecture \ref{Conj:Yoshi}:
 \begin{prop}[\cite{bib:Yos84}]
 Conjecture~\ref{Conj:Yoshi} is true, over any algebraically closed field of characteristic $0$, if $C$ is of type $\mathrm{II}$.
 \end{prop}
 Finally, we recall  that Conjecture~\ref{Conj:Yoshi} was proved in \cite[Proposition 2.7]{bib:Yos84} in the case of a nodal cubic curve, and that the group $\Aut(\Pn\backslash C)$ for this curve was studied by Wakabayashi and Yoshihara, see \cite{bib:Yos85} and \cite{bib:Wak78}.

\section{The construction}\label{Sec:ConstrDetails}
In this section, we describe precisely the construction announced in the introduction.
First we describe the triangle $\Delta$, its irreducible components and  singular points. Take three general lines of $\Pn$, that form a triangle $\Delta$, and choose the coordinates such that $\Delta$ has equation $xyz=0$. We denote by $a=(1:0:0)$, $b=(0:1:0)$, $c=(0:0:1)\in\Pn$ the singular points of $\Delta$ and by $L_{ab}$ (respectively $L_{ac}$, $L_{bc}$) the line through $a$ and $b$. In particular, $\Delta=L_{ab}\cup L_{ac}\cup L_{bc}$. 

Then, we briefly describe the two curves $\Gamma_1$ and $\Gamma_2$, in simple words. In subsection~\ref{SubSec:PointsNeighbourhoods}, we will describe these curves using the points infinitely near to $a$ and $b$.
For any $\theta\in \K^{*}$, we write $p(\theta)=(\theta:0:1)$ and denote by $\Omega_{\theta}$ the set of irreducible quartic curves of $\Pn$ that have multiplicity $3$ at $p(\theta)$, that pass through $a$ and are tangent to $L_{ab}$ and intersect $L_{bc}$ only at the point $b$.


Let $\alpha,\beta \in \K^{*}$, $\alpha\not=\beta$, then $\Gamma_1$ is one curve of $\Omega_{\alpha}$ and $\Gamma_2$ is the curve of $\Omega_{\beta}$ whose intersection with $\Gamma_1$ at the point $b$ is as large as possible. 

\subsection{The points in the neighbourhoods of $a$ and $b$}\label{SubSec:PointsNeighbourhoods}
We now describe the intersection between the curves $\Gamma_1$, $\Gamma_2$ and $\Delta$, and construct the birational morphism $\pi:X\rightarrow \Pn$ announced in Section~\ref{Sec:ConstructionIntro}.

We construct $\pi$ by a sequence of blow-ups of points that lie on the curves $\Gamma_1$, $\Gamma_2$, $\Delta$. Taking some point $x$ in a surface $S$, the blow-up $p_x:S'\rightarrow S$ gives a smooth surface $S'$. We denote by $E_x \subset S'$ the exceptional curve of $x$, which is equal to $(p_x)^{-1}(x)$. Then, $p_x$ is an isomorphism of $S'\backslash E_x$ to $S\backslash x$. It is therefore natural, for any point $y\in S\backslash x$ and any curve $C\subset S\backslash x$, to denote once again the point $p_x^{-1}(y)$ by $y$ and the curve $p_x^{-1}(C)$ by $C$. For any curve $C \subset S$ passing through $x$, the strict transform of $C$ on $S'$ will be denoted by $\widetilde{C}$. After two (or more) blow-ups, we write $\widetilde{C}=\widetilde{\widetilde{C}}$ to simplify the notation. 

\begin{figure}[ht]{
\begin{center}\includegraphics[width=10cm]{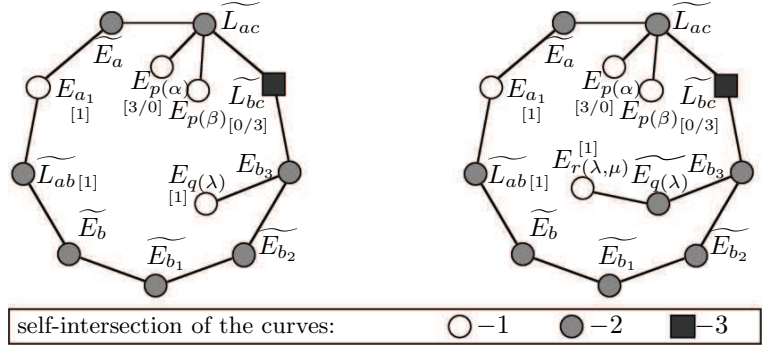}
\drawat{-35mm}{33mm}{$E_{a_1}$}%
\drawat{-30mm}{37.5mm}{$\widetilde{E_a}$}%
\drawat{-13mm}{42mm}{$\widetilde{L_{ac}}$}%
\drawat{-12mm}{32mm}{$\widetilde{L_{bc}}$}%
\drawat{-25mm}{34mm}{$E_{p(\alpha)}$}%
\drawat{-20mm}{30mm}{$E_{p(\beta)}$}%
\drawat{-11mm}{23mm}{$E_{b_3}$}%
\drawat{-8mm}{12mm}{$\widetilde{E_{b_2}}$}%
\drawat{-22.83mm}{10.77mm}{$\widetilde{E_{b_1}}$}%
\drawat{-32mm}{14.5mm}{$\widetilde{E_b}$}%
\drawat{-37.5mm}{21.5mm}{$\widetilde{L_{ab}}$}%
\drawat{-18.5mm}{21.5mm}{$\widetilde{E_{q(\lambda)}}$}%
\drawat{-29mm}{23.5mm}{${E_{r(\lambda,\mu)}}$}%
\drawat{-100mm}{1.5mm}{\small self-intersection of the curves:}%
\drawat{-39mm}{1.5mm}{$-1$}%
\drawat{-24mm}{1.5mm}{$-2$}%
\drawat{-10mm}{1.5mm}{$-3$}%
\drawat{-33mm}{30mm}{\tiny $[1]$}
\drawat{-26mm}{31mm}{\tiny $[3/0]$}
\drawat{-12mm}{28.5mm}{\tiny $[0/3]$}
\drawat{-25.5mm}{25.5mm}{\tiny $[1]$}
\drawat{-32mm}{21mm}{\tiny $[1]$}
\drawat{-95mm}{33mm}{$E_{a_1}$}%
\drawat{-90mm}{37.5mm}{$\widetilde{E_a}$}%
\drawat{-73mm}{42mm}{$\widetilde{L_{ac}}$}%
\drawat{-72mm}{32mm}{$\widetilde{L_{bc}}$}%
\drawat{-85mm}{34mm}{$E_{p(\alpha)}$}%
\drawat{-80mm}{30mm}{$E_{p(\beta)}$}%
\drawat{-71mm}{23mm}{$E_{b_3}$}%
\drawat{-68mm}{12mm}{$\widetilde{E_{b_2}}$}%
\drawat{-82.83mm}{10.77mm}{$\widetilde{E_{b_1}}$}%
\drawat{-92mm}{14.5mm}{$\widetilde{E_b}$}%
\drawat{-97.5mm}{21.5mm}{$\widetilde{L_{ab}}$}%
\drawat{-80mm}{21.5mm}{$E_{q(\lambda)}$}%
\drawat{-93mm}{30mm}{\tiny $[1]$}
\drawat{-86mm}{31mm}{\tiny $[3/0]$}
\drawat{-72mm}{28.5mm}{\tiny $[0/3]$}
\drawat{-80mm}{19mm}{\tiny $[1]$}
\drawat{-92mm}{21mm}{\tiny $[1]$}
\end{center}

} \caption{\small The situation on the surface $X$. Two curves are connected by an edge if their intersection is positive (and here equal to $1$). The positive intersections with $\widetilde{\Gamma_1}$ and $\widetilde{\Gamma_2}$ are in square brackets.\label{SituationX}}\end{figure}

Our aim is to obtain the configuration of curves of Figure \ref{SituationX} on $X$. For this, we will blow-up the points $p(\alpha)$, $p(\beta)$, and points in the neighbourhoods of $a$ and $b$.

\begin{figure}[ht]{
\begin{center}\includegraphics[width=10cm]{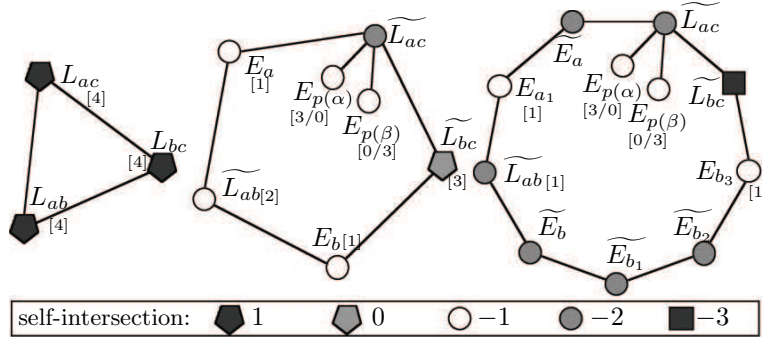}
\drawat{-94.5mm}{33.68mm}{$L_{ac}$}%
\drawat{-82.5mm}{25mm}{$L_{bc}$}%
\drawat{-98.6mm}{17.5mm}{$L_{ab}$}%
\drawat{-70mm}{35mm}{$E_a$}%
\drawat{-51mm}{39mm}{$\widetilde{L_{ac}}$}%
\drawat{-44mm}{25.5mm}{$\widetilde{L_{bc}}$}%
\drawat{-61mm}{12mm}{$E_b$}%
\drawat{-73mm}{18.71mm}{$\widetilde{L_{ab}}$}%
\drawat{-63.5mm}{31.5mm}{$E_{p(\alpha)}$}%
\drawat{-57mm}{27mm}{$E_{p(\beta)}$}%
\drawat{-34mm}{32mm}{$E_{a_1}$}%
\drawat{-29mm}{36.5mm}{$\widetilde{E_a}$}%
\drawat{-12mm}{41mm}{$\widetilde{L_{ac}}$}%
\drawat{-11mm}{31mm}{$\widetilde{L_{bc}}$}%
\drawat{-25mm}{32.5mm}{$E_{p(\alpha)}$}%
\drawat{-19mm}{28.5mm}{$E_{p(\beta)}$}%
\drawat{-10mm}{21mm}{$E_{b_3}$}%
\drawat{-13mm}{13.5mm}{$\widetilde{E_{b_2}}$}%
\drawat{-21.83mm}{9.77mm}{$\widetilde{E_{b_1}}$}%
\drawat{-31mm}{13.5mm}{$\widetilde{E_b}$}%
\drawat{-35.5mm}{20.5mm}{$\widetilde{L_{ab}}$}%
\drawat{-100mm}{1.5mm}{\small self-intersection:}%
\drawat{-69mm}{1.5mm}{$1$}%
\drawat{-53mm}{1.5mm}{$0$}%
\drawat{-39mm}{1.5mm}{$-1$}%
\drawat{-24mm}{1.5mm}{$-2$}%
\drawat{-10mm}{1.5mm}{$-3$}%
\drawat{-91mm}{31mm}{\tiny $[4]$}
\drawat{-96mm}{14.2mm}{\tiny $[4]$}
\drawat{-85.5mm}{22.5mm}{\tiny $[4]$}
\drawat{-69.23mm}{32.85mm}{\tiny $[1]$}
\drawat{-68mm}{18.09mm}{\tiny $[2]$}
\drawat{-57mm}{12mm}{\tiny $[1]$}
\drawat{-43mm}{20mm}{\tiny $[3]$}
\drawat{-55mm}{24mm}{\tiny $[0/3]$}
\drawat{-64mm}{28.5mm}{\tiny $[3/0]$}
\drawat{-33mm}{29mm}{\tiny $[1]$}
\drawat{-25mm}{29.5mm}{\tiny $[3/0]$}
\drawat{-19mm}{25.5mm}{\tiny $[0/3]$}
\drawat{-3mm}{18.5mm}{\tiny $[1]$}
\drawat{-30mm}{19.96mm}{\tiny $[1]$}
\end{center}
} \caption{\small The configuration of the special curves on the surface  $X'$\label{ConfigXprime}. Two curves are connected by an edge if their intersection is positive (and here equal to $1$). The positive intersections with $\widetilde{\Gamma_1}$ and $\widetilde{\Gamma_2}$ are in square brackets.}\end{figure}

Denote by $a_1$ the point in the first neighbourhood of $a$ that belongs to the (strict transform of the) line $L_{ab}$, and by $b_1$ the point in the first neighbourhood of $b$ that belongs to the line $L_{bc}$. For $i=2,3$, we call $b_i$ the point in the first neighbourhood of  $b_{i-1}$ (and thus in the $i$-th neighbourhood of $b$) that belongs to the line $L_{bc}$. We denote by $\pi':X'\rightarrow \Pn$ the blow-up of the points $a$, $a_1$, $b$, $b_1$, $b_2$, $b_3$, $p(\alpha)$ and $p(\beta)$. The configuration of the curves on $X'$ and the decomposition of $\pi'$  are described in Figure~\ref{ConfigXprime}.

On the surface $X'$, (the strict pull-back of) any curve of $\Omega_{\alpha}$ has self-intersection~$1$, and its intersection with $E_{p(\alpha)}$, ${E_{a_1}}$, $E_{b_3}$ and $\widetilde{L_{ab}}$  is respectively $3$, $1$, $1$ and $1$; furthermore no other curve of Figure~\ref{ConfigXprime} intersects any curve of $\Omega_{\alpha}$. The situation for the curves of $\Omega_{\beta}$ is similar, after exchanging the roles of $E_{p(\alpha)}$ and $E_{p(\beta)}$.

Since $E_{b_3}\cong \mathbb{P}^1$, the points of $E_{b_3}$ that do not lie on $\widetilde{L_{bc}}$ or $\widetilde{E_{b_2}}$ are parametrised by $\K^{*}$.
Explicitly, the morphism $\pi':X'\rightarrow \Pn$ is given locally by $(x,y)\mapsto (xy^4:1:y)$, and in these coordinates, we define for any $\theta \in \K^{*}$ the point $q(\theta)\in E_{b_3}\subset X'$ that corresponds to $(\theta,0)$. Any curve of $\Omega_{\alpha}$ (respectively of $\Omega_{\beta}$) passes through $q(\theta)$, for some $\theta \in \K^{*}$.

We assume that both $\Gamma_1$ and $\Gamma_2$ pass through the same point $q(\lambda)\in X'$, which is consistent with the fact that $\Gamma_1$ and $\Gamma_2$ have their maximum intersection at $b$. Blowing-up $q(\lambda)$, the exceptional curve $E_{q(\lambda)}$  intersects $\widetilde{E_{b_3}}$  in one point, through which no curve of $\Omega_{\alpha}$ or $\Omega_{\beta}$ passes. The remaining points of $E_{q(\lambda)}$ are parametrised by $\K$.
Using the same coordinates as above, the blow-up of $q(\lambda)=(\lambda,0)$ may be viewed as $(x,y)\mapsto (xy+\lambda,y)$, and the parametrisation associates to $\mu\in \K$ the point $r(\lambda,\mu)$, equal to $(\mu,0)$.

\begin{lemm}\label{Lem:UniqueCurve}
For any pair $(\lambda,\mu)\in \K^{*}\times \K$, there exists a unique curve in $\Omega_{\alpha}$, that passes through $q(\lambda)$ and $r(\lambda,\mu)$. The same is true for $\Omega_{\beta}$. The equations of the two curves are
\begin{equation}\begin{array}{l}
\lambda^2z(\alpha z-x)^3+\alpha^2xy^2 (\mu (\alpha z-x)-\alpha\lambda y),\\
\lambda^2z(\beta z-x)^3+\beta^2xy^2 (\mu (\beta z-x)-\beta\lambda y).\end{array}\end{equation}
\end{lemm}
\begin{proof}
This follows from a straightforward calculation, using the description  of the blow-up in coordinates given above.
\end{proof}

From now on, we fix $(\lambda,\mu)\in \K^{*}\times \K$, and denote by $\Gamma_1\subset \Omega_{\alpha}$ and $\Gamma_2\in\Omega_{\beta}$ the two curves yielded by Lemma~\ref{Lem:UniqueCurve}. Blowing-up the point $q(\lambda)$ on $X'$ and then the point $r(\lambda,\mu)$, we obtain the birational morphism $\pi:X\rightarrow \Pn$ announced in the introduction. The situation on the blow-up of $X'$ at $q(\lambda)$ and on the surface $X$ is described in Figure~\ref{SituationX}.

On the surface $X$, let $R$ be the reducible curve which is the union of the $9$ curves of self-intersection $\leq -2$ of Figure~\ref{SituationX} (the curves in grey).

\begin{prop}\label{Prp:IsoCompl}
Fix some $i\in\{1,2\}$. There exists a birational morphism $\eta_i:X\rightarrow \Pn$ that collapses the curves $R\cup \widetilde{\Gamma_i}$; it starts by collapsing $\widetilde{\Gamma_i}$ and then collapses the images of respectively $\widetilde{L_{ab}}$, $\widetilde{E_b}$, $\widetilde{E_{b_1}}$, $\widetilde{E_{b_2}}$, $\widetilde{E_{b_3}}$, $\widetilde{E_{q(\lambda)}}$, $\widetilde{L_{bc}}$, $\widetilde{L_{ac}}$, $\widetilde{E_a}$.

Then,  $\eta_i(\widetilde{\Gamma_{3-i}})$ is a curve of $\Pn$ of  degree $39$, irreducible over the algebraic closure of $\K$,  which has exactly one singular point. The morphism $\eta_i$ is a minimal resolution of this curve, and is a $(-1)$-tower resolution of it (see Definition~\ref{Defi:TowerResolution}). 
\end{prop}
\begin{proof}
The curve $\widetilde{\Gamma_i}$ is a $(-1)$-curve (a smooth rational curve of self-intersection $-1$, irreducible over the algebraic closure of $\K$). We may therefore collapse it and obtain a birational morphism $X\rightarrow Y$ where $Y$ is smooth and projective. On $Y$, the image of $\widetilde{L_{ab}}$ is a $(-1)$-curve so we may collapse it. Continuing with the images of $\widetilde{E_b}$, $\widetilde{E_{b_1}}$, ..., $\widetilde{E_a}$ we obtain a birational morphism $\eta_i:X\rightarrow Z$ for some smooth rational projective surface $Z$ (see Figure \ref{descente}).

\begin{figure}[ht]{
\begin{center}\includegraphics[width=10cm]{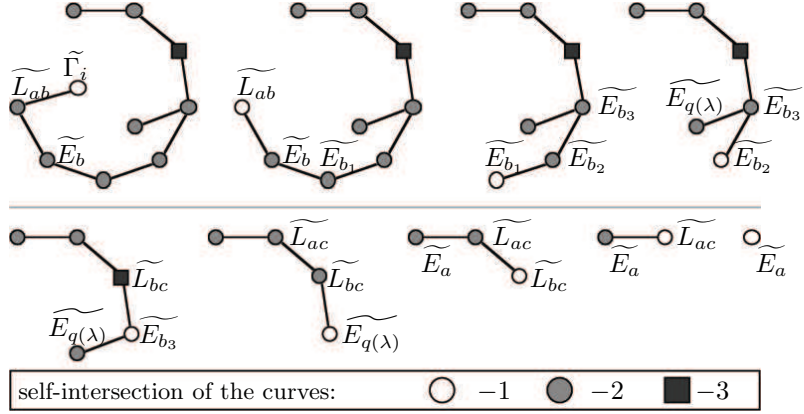}
\drawat{-100mm}{1.5mm}{\small self-intersection of the curves:}%
\drawat{-39mm}{1.5mm}{$-1$}%
\drawat{-24mm}{1.5mm}{$-2$}%
\drawat{-10mm}{1.5mm}{$-3$}%
\drawat{-95mm}{33mm}{$\widetilde{E_b}$}%
\drawat{-94mm}{44.5mm}{$\widetilde{\Gamma_i}$}%
\drawat{-101mm}{42mm}{$\widetilde{L_{ab}}$}%
\drawat{-71mm}{42mm}{$\widetilde{L_{ab}}$}%
\drawat{-65mm}{33mm}{$\widetilde{E_b}$}%
\drawat{-60mm}{32.5mm}{$\widetilde{E_{b_1}}$}%
\drawat{-38mm}{32.5mm}{$\widetilde{E_{b_1}}$}%
\drawat{-27mm}{32.5mm}{$\widetilde{E_{b_2}}$}%
\drawat{-23mm}{39.5mm}{$\widetilde{E_{b_3}}$}%
\drawat{-5mm}{32.5mm}{$\widetilde{E_{b_2}}$}%
\drawat{-1mm}{39.5mm}{$\widetilde{E_{b_3}}$}%
\drawat{-14mm}{40.5mm}{$\widetilde{E_{q(\lambda)}}$}%
\drawat{-84mm}{9mm}{$\widetilde{E_{b_3}}$}%
\drawat{-96mm}{10mm}{$\widetilde{E_{q(\lambda)}}$}%
\drawat{-85mm}{16.5mm}{$\widetilde{L_{bc}}$}%
\drawat{-59mm}{16.5mm}{$\widetilde{L_{bc}}$}%
\drawat{-57mm}{9mm}{$\widetilde{E_{q(\lambda)}}$}%
\drawat{-64mm}{22mm}{$\widetilde{L_{ac}}$}%
\drawat{-37mm}{22mm}{$\widetilde{L_{ac}}$}%
\drawat{-32mm}{16.5mm}{$\widetilde{L_{bc}}$}%
\drawat{-46.5mm}{18.5mm}{$\widetilde{E_a}$}%
\drawat{-21.5mm}{18.5mm}{$\widetilde{E_a}$}%
\drawat{-12.5mm}{22mm}{$\widetilde{L_{ac}}$}%
\drawat{-2mm}{18.5mm}{$\widetilde{E_a}$}%
\end{center}

} \caption{\small The decomposition of the birational morphism $X\rightarrow Z$ and the image of $R\cup \widetilde{\Gamma_i}$ on each surface.\label{descente}}\end{figure}
 Since $X$ was obtained by blowing-up $10$ points from $\Pn$ and $\eta_i$ collapses $10$ irreducible curves, we have $(K_Z)^2=(K_{\Pn})^2=9$, so $Z\cong \Pn$.

Write $j=3-i$. Since $\widetilde{\Gamma_{j}}$ is not collapsed by $\eta_i$, the image $\eta_i(\widetilde{\Gamma_{j}})$ is a curve. Its irreducibility follows from that of $\widetilde{\Gamma_{j}}$. Its degree can be calculated by computing its self-intersection after each of the $10$ blow-downs. 
Since $R\cup \widetilde{\Gamma_i}$ is connected, its image by $\eta_i$ is a single point. The curve $\widetilde{\Gamma_j}$ is smooth and intersects $\widetilde{\Gamma_i}$ in more than one point, hence $\eta_i$ is a minimal resolution of $\eta_i(\widetilde{\Gamma_{j}})$ and this curve has a unique singular point. Furthermore $\eta_i$ is a tower resolution, as it collapses only one curve of self-intersection $-1$.\end{proof}

\section{Comparison of the curves $\eta_1(\widetilde{\Gamma_{2}})$ and $\eta_2(\widetilde{\Gamma_{1}})$}\label{Sec:Comparison}
Proposition~\ref{Prp:IsoCompl} shows that for any choice of $\alpha, \beta,\lambda\in \K^{*}, \alpha\not=\beta, \mu\in\K$, the complements of the two curves $\eta_1(\widetilde{\Gamma_{2}})$ and $\eta_1(\widetilde{\Gamma_{2}})$ are isomorphic. In this section, we distinguish the differences between the two curves. 
\begin{prop}\label{Prp:ProjEquiv}
The following are equivalent:
\begin{enumerate}
\item
there exists an automorphism $\psi$ of $\Pn$ that sends $\eta_1(\widetilde{\Gamma_{2}})$ on $\eta_2(\widetilde{\Gamma_{1}})$;
\item
there exists an automorphism $\psi'$ of $X$ that leaves invariant every irreducible component of $R$ and exchanges $\widetilde{\Gamma_1}$ and $\widetilde{\Gamma_2}$;
\item
there exists an automorphism $\psi''$ of $\Pn$ that fixes $a$, $b$ and $c$ and permutes $\Gamma_1$ and $\Gamma_2$;
\item
$\mu=0$ and $\alpha+\beta=0$.
\end{enumerate}
\end{prop}
\begin{proof}
Let us keep Diagram~\ref{Eq:DiagrammFond} in mind.
The fact that $\eta_1$ (respectively $\eta_2$) is a minimal resolution of $\eta_{1}(\widetilde{\Gamma_{2}})$ (respectively of $\eta_{2}(\widetilde{\Gamma_1})$) and the assumptions made on the automorphisms above imply that  $\psi'$ may be constructed starting from $\psi$, as $\psi'=\eta_1^{-1}\psi\eta_2$. Similarly, the existence of $\psi'$ implies that of $\psi$ and $\psi''$, constructed as $\psi=\eta_1\psi'\eta_2^{-1}$ and $\psi''=\pi\psi'\pi^{-1}$. Finally, if $\psi''$ exists, then $\psi'=\pi^{-1}\psi''\pi$ exists.
\begin{center}
$\xymatrix@R=0.5mm@C=2cm{& &  \Pn\\
\Pn\ar@(lu,ru)^{\psi''}&X\ar@(lu,ru)^{\psi'}\ar[l]_{\pi}\ar@/_/[rd]^{\eta_2}\ar@/^/[ru]^{\eta_1}&\\
&&\Pn\ar@/^1pc/@{-->}[uu]_{\varphi}
\ar@/_1pc/[uu]_{\psi} }$\end{center}
It remains to prove that assertions $(3)$ and $(4)$ are equivalent. If  $\psi''$ exists, then it is of the form $(x:y:z)\mapsto (x:\xi y: \theta z)$, for some $\xi,\theta\in\K^{*}$. Since $\psi''$ exchanges the curves $\Gamma_1$ and $\Gamma_2$, it exchanges the points $p(\alpha)=(\alpha:0:1)$ and $p(\beta)=(\beta:0:1)$, which implies that $\alpha+\beta=0$ and $\theta=-1$. Using the explicit equations of $\Gamma_1$ and $\Gamma_2$, we find directly that $\mu=0$. 
Conversely, if $\mu=0$ and $\alpha+\beta=0$ the automorphism $(x:y:z)\mapsto (x:y:-z)$ exchanges $\Gamma_1$ and $\Gamma_2$.
\end{proof}
Propositions~\ref{Prp:IsoCompl} and \ref{Prp:ProjEquiv} yield counterexamples to Conjecture~\ref{Conj:Yoshi}, for any field~$\K$ that has more than two elements.
Now, we study more intrinsically the curves $\eta_1(\widetilde{\Gamma_2})$, $\eta_2(\widetilde{\Gamma_1})$, without taking in account the plane embedding.
\begin{prop}\label{Prp:NonIsomorphic}
Neither $\eta_1(\widetilde{\Gamma_2})$ nor $\eta_2(\widetilde{\Gamma_1})$ is a curve of type $\mathrm{I}$.

If $\mu=0$, the curves  $\eta_1(\widetilde{\Gamma_2})$ and $\eta_2(\widetilde{\Gamma_1})$ are isomorphic.

For any field~$\K$ with more than two elements there exist values of $\alpha,\beta,\mu$ for which the curves  $\eta_1(\widetilde{\Gamma_2})$ and $\eta_2(\widetilde{\Gamma_1})$ are not isomorphic.
\end{prop}
\begin{proof}
Denote by $p_1$ (respectively $p_2$) the morphism $\widetilde{\Gamma_1}\rightarrow \eta_2(\widetilde{\Gamma_1})$ (respectively $\widetilde{\Gamma_2}\rightarrow \eta_1(\widetilde{\Gamma_2})$) obtained by restriction of $\eta_2$ (respectively $\eta_1$).
The singular curves $\eta_1(\widetilde{\Gamma_2})$ and $\eta_2(\widetilde{\Gamma_1})$ are isomorphic if and only if there is an isomorphism $\rho:\widetilde{\Gamma_1}\rightarrow\widetilde{\Gamma_2}$ that is compatible with $p_1$ and $p_2$. Furthermore the singular curves are of type $I$ if and only if the morphisms $p_i$ are injective.

The ramified form of the morphism $p_1$ consists of the point $\widetilde{L_{ab}}\cap \widetilde{\Gamma_1}$ and the form of degree $8$ on $\widetilde{\Gamma_1}\cong\mathbb{P}^1$ obtained by intersecting $\widetilde{\Gamma_1}$ with $\widetilde{\Gamma_2}$. Taking some coordinates $(u:v)$ on $\widetilde{\Gamma_1}\cong\mathbb{P}^1$, the morphism $\widetilde{\Gamma_1}\rightarrow \Gamma_1\subset\Pn$ obtained by restriction of $\pi$ is the following:
\[(u:v)\mapsto \left(v^4\lambda^5\alpha:(u+\mu v)(\alpha(u+\mu v)^2u-\lambda^4v^3):v(u+\mu v)^2\lambda u\alpha\right).\]
The point $(1:0)$ is sent on $b$, the point $(\mu:-1)$ is sent on $a$ and the point $(0:1)$ corresponds to $\widetilde{\Gamma_1}\cap \widetilde{L_{ab}}$.
Replacing the parametrisation in the equation of $\Gamma_1$ we find $0$, and replacing it in the equation of $\Gamma_2$, we find 
\begin{center}$-\lambda^{10}\alpha(\alpha-\beta)(u+\mu v)^2 v^7(\sum_{i=0}^7 c_i\cdot  u^iv^{7-i})$,\end{center}
where
 $c_0,...,c_7$ are as follows:
\begin{center}$\begin{array}{rclrcl}
c_0&=&3\alpha^2\beta^2 & c_4&=&-\alpha\beta\mu(8\lambda^4\beta-7\alpha\beta\mu^3+6\lambda^4\alpha)\\
c_1&=&13\alpha^2\beta^2\mu& c_5&=& -\alpha\beta\mu^2(3\lambda^4\alpha-\alpha\beta\mu^3+7\lambda^4\beta)\\
c_2&=&22\alpha^2\beta^2\mu^2 & c_6&=&\lambda^4(\lambda^4(\alpha\beta+\alpha^2+\beta^2)-2\alpha\beta^2\mu^3 \\
c_3&=&-3\alpha\beta(\lambda^4(\alpha+\beta)-6\alpha\beta\mu^3)&c_7&=&\lambda^8\beta^2\mu.
\end{array}$\end{center}
The intersection number of $\Gamma_1$ and $\Gamma_2$ is $16$; the intersections at $a$ and $a_1$ correspond to the factor $(u+\mu v)^2$ and the intersections at $b, b_1,b_2,b_3, q(\lambda),r(\lambda,\mu)$ correspond to $v^6$. Thus, the form of degree $8$ on $\widetilde{\Gamma_1}$ corresponding to the intersection of $\widetilde{\Gamma_1}$ and $\widetilde{\Gamma_2}$ is $F_1=v\sum_{i=0}^7 c_i\cdot  u^iv^{7-i}$. Since $F_1$ vanishes at the point $(1:0)$ and $(0:1)$ corresponds to $\widetilde{\Gamma_1}\cap \widetilde{L_{ab}}$, the map $p_1$ is not injective and $\eta_2(\widetilde{\Gamma_1})$ is not of type $\mathrm{I}$.

For $p_2:\widetilde{\Gamma_2}\rightarrow \eta_1(\widetilde{\Gamma_2})$, the situation is similar. We find a form $F_2$, that is equal to $F_1$, after exchanging $\alpha$ and $\beta$. As above, we see that $\eta_1(\widetilde{\Gamma_2})$ is not of type $\mathrm{I}$. Finally, the two singular curves are isomorphic if and only if there exists an isomorphism of $\mathbb{P}^1$ that fixes $(0:1)$ and sends $F_1$ on $F_2$ (we say in this case that $F_1$ and $F_2$ are equivalent). If $\mu=0$, the identity suits, since each $c_i$ becomes symmetric with respect to $\alpha$ and $\beta$. If $\mu\not=0$, this is not the case. If $\car(\K)\not=2$, choosing $\alpha=1,\beta=2,\lambda=\mu=1$, we can compute that $F_1$ and $F_2$ are not equivalent. If $\car(\K)=2$, there is considerable simplification of the terms, and we find that if $\lambda=\mu=1$, $\alpha\not=\beta$, then $F_1$ and $F_2$ are not equivalent.\end{proof}

The proof of Theorems~\ref{Thm:CountExConj} and \ref{Thm:SecondConj} now follows directly from Propositions~\ref{Prp:IsoCompl}, \ref{Prp:ProjEquiv} and \ref{Prp:NonIsomorphic}.

\end{document}